\documentclass[11pt,reqno]{article}

\usepackage{amsmath,amsthm,amsfonts,amssymb,latexsym,mathrsfs,color}

\newtheorem{theorem}{Theorem}[section]
\newtheorem{lemma}[theorem]{Lemma}
\newtheorem{corollary}[theorem]{Corollary}
\newtheorem{proposition}[theorem]{Proposition}
\newtheorem{problem}[theorem]{Problem}
\newtheorem{conjecture}[theorem]{Conjecture}
\newtheorem{question}[theorem]{Question}

\theoremstyle{definition}
\newtheorem{definition}[theorem]{Definition}
\newtheorem{example}[theorem]{Example}

\theoremstyle{remark}
\newtheorem{remark}[theorem]{Remark}

\numberwithin{equation}{section}

\newcommand{\la}{\lambda}
\newcommand{\mo}{\mathcal{O}}
\newcommand{\msn}{{\mathcal S}_n}
\newcommand{\exc}{{\rm exc\,}}
\newcommand{\bn}{{\rm\bf n}}
\newcommand{\sgn}{{\rm sgn\,}}
\newcommand{\rz}{{\rm RZ}}
\newcommand{\pf}{{\rm PF}}
\newcommand{\sep}{\preceq}
\newcommand{\seps}{\prec}

\title{A unified approach to \\ polynomial sequences with only real zeros
\thanks{Partially supported by NSF of China 10471016.}}
\author{Lily L. Liu,\quad Yi Wang
\footnote{Corresponding author.
\newline\hspace*{5mm}
   {\it Email addresses:}\quad lliulily@yahoo.com.cn (L.L. Liu), wangyi@dlut.edu.cn (Y. Wang).}}
\date{\footnotesize Department of Applied Mathematics,
         Dalian University of Technology,
         Dalian 116024, P. R. China}

\begin{document}

\maketitle

\begin{abstract}
We give new sufficient conditions for a sequence of polynomials to
have only real zeros based on the method of interlacing zeros. As
applications we derive several well-known facts, including the
reality of zeros of orthogonal polynomials, matching polynomials,
Narayana polynomials and Eulerian polynomials. We also settle
certain conjectures of Stahl on genus polynomials by proving them
for certain classes of graphs, while showing that they are false in
general.
\end{abstract}
\section{Introduction}
\hspace*{\parindent}
Polynomials with only real zeros arise often in combinatorics and
other branches of mathematics (see
\cite{BraTAMS,Bre89,Bre94,Hag00,Pit97,Sta89,Sta00,Wag91,Wag92,WYjcta05}).
Since our interest is combinatorial, we are mainly concerned with
polynomials whose coefficients are positive or alternating in sign.
The real zeros of such polynomials are either all negative or all
positive. A polynomial with coefficients alternating in sign can be
converted into a polynomial with positive coefficients. So we may
concentrate our attention on polynomials with positive coefficients.
We were led to the study of such polynomials because of their
implications on unimodality and log-concavity.

Let $a_0,a_1,\ldots,a_n$ be a sequence of positive real numbers. The
sequence is said to be {\it unimodal} if there exists an index $0\le
m\le n$, called the {\it mode} of the sequence, such that
$a_0\le\cdots\le a_{m-1}\le a_m\ge a_{m+1}\ge\cdots\ge a_n$. The
sequence is said to be {\it log-concave} if $a_{i-1}a_{i+1}\le
a_i^2$ for $i=1,\ldots,n-1$. Clearly, log-concavity implies
unimodality. Unimodal and log-concave sequences occur naturally in
combinatorics, analysis, algebra, geometry, probability and
statistics. The reader is referred to Stanley~\cite{Sta89} and
Brenti~\cite{Bre94} for surveys and
\cite{Bre95,Bre96,wyjcta02,wyeujc02,wylaa03,WYeujc05,WYjcta05,WYlcp}
for recent progress on this subject.

One classical approach to unimodality and log-concavity of a finite sequence
is to use Newton's inequality:
if the polynomial $\sum_{i=0}^{n}a_ix^i$ with positive coefficients has only real zeros,
then
$$a_i^2\ge a_{i-1}a_{i+1}\left(1+\frac{1}{i}\right)\left(1+\frac{1}{n-i}\right)$$
for $1\le i\le n-1$, and the sequence $a_0,a_1,\ldots,a_n$ is therefore unimodal and log-concave
(see Hardy, Littlewood and P\'olya~\cite[p. 104]{HLP52}).
Such a sequence of positive numbers whose generating function has only real zeros
is called a {\it P\'olya frequency sequence} in the theory of total positivity.
See Karlin~\cite{Kar68} for a standard reference on total positivity
and Brenti~\cite{Bre89,Bre95,Bre96}
for applications of total positivity to unimodality and log-concavity problems.
It often occurs that unimodality of a sequence is known,
yet to determine the exact number and location of modes
is a much more difficult task.
The case for P\'olya frequency sequences is somewhat different.
Darroch~\cite{Dar64} showed that
if the polynomial $P(x)=\sum_{i=0}^{n}a_ix^i$ with positive coefficients has only real zeros,
then the unimodal sequence $a_0,a_1,\ldots,a_n$ has at most two modes
and each mode $m$ satisfies
$$\left\lfloor\frac{P'(1)}{P(1)}\right\rfloor\le m\le \left\lceil\frac{P'(1)}{P(1)}\right\rceil.$$

Our main concern here is sequences of polynomials with only real zeros.
Such polynomial sequences occurring in combinatorics often satisfy certain recurrence relations.
For example,
the sequence of classical orthogonal polynomials $p_n(x)$ satisfies a three-term recurrence relation
\begin{equation}\label{op-rr}
p_n(x)=(a_nx+b_n)p_{n-1}(x)-c_np_{n-2}(x)
\end{equation}
with $p_{-1}(x)=0$ and $p_0(x)=1$,
where $a_n,c_n>0$ and $b_n\in\mathbb{R}$
(see \cite[Theorem 3.2.1]{Sze75}).
A standard result in the theory of orthogonal polynomials
is that for each $n\ge 1$,
the zeros of $p_n(x)$ are real, simple, and separate those of $p_{n+1}(x)$.
In this paper, we develop methods to provide a unified approach to the reality of zeros of
polynomial sequences satisfying certain recurrence relations.

Following Wagner~\cite{Wag92},
a real polynomial is said to be {\it standard}
if either it is identically zero or its leading coefficient is positive.
Let $\rz$ denote the set of real polynomials with only real zeros.
In particular,
let $\pf$ consist of those polynomials in $\rz$ whose coefficients are nonnegative.
In other words,
$\pf$ is the set of polynomials whose coefficients form a P\'olya frequency sequence.
Clearly, all zeros of a polynomial in $\pf$ are real and nonpositive.
For convenience, let $0\in\pf$.

Suppose that $f,g\in\rz$. Let $\{r_i\}$ and $\{s_j\}$ be all zeros
of $f$ and $g$ in nonincreasing order respectively. We say that $g$
{\it alternates left of} $f$ ($g$ alternates $f$ for short) if $\deg
f=\deg g=n$ and
\begin{equation}\label{alt-def}
s_n\le r_n\le s_{n-1}\le\cdots\le s_2\le r_2\le s_1\le r_1.
\end{equation}
We say that $g$ {\it interlaces} $f$
if $\deg f=\deg g+1=n$ and
\begin{equation}\label{int-def}
r_n\le s_{n-1}\le\cdots\le s_2\le r_2\le s_1\le r_1.
\end{equation}
Let $g\sep f$ denote ``either $g$ alternates $f$ or $g$ interlaces $f$''.
If no equality sign occurs in (\ref{alt-def}) (resp. (\ref{int-def})),
then we say that $g$ strictly alternates $f$ (resp. $g$ strictly interlaces $f$).
Let $g\seps f$ denote ``either $g$ strictly alternates $f$ or $g$ strictly interlaces $f$''.
For notational convenience, let $a\sep bx+c$ for any real constants $a,b,c$
and $f\sep 0, 0\sep f$ for any real polynomial $f$ with only real zeros.

Let $\{P_n(x)\}_{n\ge 0}$ be a sequence of standard polynomials.
We say that $\{P_n(x)\}$ is a {\it Sturm sequence}
if $\deg P_n=n$ and $P_{n-1}(r)P_{n+1}(r)<0$ whenever $P_n(r)=0$ and $n\ge 1$.
It is well known that $\{P_n(x)\}$ is a Sturm sequence if and only if,
for $n\ge 1$, $P_n\in\rz$ and $P_{n}$ strictly interlaces $P_{n+1}$.
We say that $\{P_n(x)\}$ is a {\it generalized Sturm sequence} if
$P_n\in\rz$ and $P_0\sep P_1\sep\cdots\sep P_{n-1}\sep P_n\sep\cdots$.
For example,
if $P$ is a standard polynomial with only real zeros and $\deg P=n$,
then $P^{(n)},P^{(n-1)},\ldots,P',P$ form a generalized Sturm sequence
by Rolle's theorem.
In particular, if the zeros of $P$ are distinct,
then $P^{(n)},P^{(n-1)},\ldots,P',P$ form a Sturm sequence.

There are various methods for showing that polynomials have only real zeros.
For example, Wang and Yeh~\cite[Theorem 1]{WYjcta05}
established the following result
which has been proved to be an extremely useful tool.
In fact, it provides a unified approach to
unimodality and log-concavity of many well-known sequences in combinatorics.
See \cite{WYjcta05} for details.
\begin{theorem}\label{WY-thm}
Let $f$ and $g$ be real polynomials whose leading coefficients have the same sign.
Suppose that $f,g\in\rz$ and $g\sep f$.
If $ad\le bc$,
then $(ax+b)f(x)+(cx+d)g(x)\in\rz$.
\end{theorem}

An immediate consequence of Theorem~\ref{WY-thm} is the following corollary
(\cite[Corollary 1]{WYjcta05}).
\begin{corollary}\label{WY-cor}
Suppose that $f,g\in\pf$ and $g$ interlaces $f$.
If $ad\ge bc$, then $(ax+b)f(x)+x(cx+d)g(x)\in\rz$.
\end{corollary}

Theorem~\ref{WY-thm} and Corollary~\ref{WY-cor} provide the inductive basis
for showing the reality of zeros of polynomial sequences $\{P_n(x)\}$
satisfying certain recurrence relations $P_n(x)=a_n(x)P_{n-1}(x)+b_n(x)P'_{n-1}(x)$.
For example,
let $S(n,k)$ be the Stirling number of the second kind and
let $S_n(x)=\sum_{k=0}^{n}S(n,k)x^k$ denote the associated generating function.
It is well known that
$S(n,k)=S(n-1,k-1)+kS(n-1,k)$, $S(0,0)=1$ and $S(n,0)=0$ for $n\ge 1$,
which is equivalent to
\begin{equation}\label{snx-rr}
S_n(x)=xS_{n-1}(x)+xS'_{n-1}(x)
\end{equation}
with $S_0(x)=1$.
From Corollary~\ref{WY-cor} (or Theorem~\ref{WY-thm}) and by induction,
it follows immediately that $S_n(x)$ has only real zeros for each $n\ge 1$,
a result originally due to Harper~\cite{Har67}.
Similarly,
the classical Eulerian polynomial $A_n(x)$,
which satisfies the recurrence relation
\begin{equation}\label{anq-rr}
A_n(x)=nxA_{n-1}(x)+x(1-x)A'_{n-1}(x)
\end{equation}
with $A_0(x)=1$
(see Comtet~\cite[Exercise VII. 3]{Com74} for instance),
has only real zeros for $n\ge 1$.

However, Theorem~\ref{WY-thm} and Corollary~\ref{WY-cor} cannot be used to prove the fact
that $\{S_n(x)\}$ and $\{A_n(x)\}$ form generalized Sturm sequences respectively.
Neither are applicable to polynomial sequences $\{P_n(x)\}$
satisfying recurrence relations similar to
$P_n(x)=a_n(x)P_{n-1}(x)+b_n(x)P_{n-2}(x)$,
nor even to the orthogonal polynomials.
So a natural problem arises:
Given $F(x)=a(x)f(x)+b(x)g(x)$,
under what conditions $g\sep f$ necessarily implies $f\sep F$?
This leads to the following two more general problems.
\begin{problem}\label{prob-Ffgk}
Let $F(x)=a(x)f(x)+\sum_{j=1}^{k}b_j(x)g_j(x)$.
Suppose that $f,g_j\in\rz$ and $g_j\sep f$ for all $j$.
Under what conditions $F\in\rz$ and $f\sep F$?
\end{problem}
\begin{problem}\label{prob-FGfg}
Let $F(x)=a(x)f(x)+b(x)g(x)$ and $G(x)=c(x)f(x)+d(x)g(x)$.
Suppose that $f,g\in\rz$ and $g\sep f$.
Under what conditions $F,G\in\rz$ and $G\sep F$?
\end{problem}

The organization of this paper is as follows.
In Section 2, we present various results
concerning Problem~\ref{prob-Ffgk} and \ref{prob-FGfg}.
Among other things, we give a simple proof of Theorem~\ref{WY-thm}.
In Section 3, we apply these results to derive several well-known facts
and solve certain open problems and conjectures
in a unified manner.
\section{Main results}
\hspace*{\parindent}
Let $\sgn$ denote the sign function defined on $\mathbb{R}$ by
$$\sgn(x)=\begin{cases}
+1 & \text{if $x>0$},\\ 0 & \text{if $x=0$},\\ -1 & \text{if $x<0$}.
\end{cases}$$
Let $f(x)$ be a real function.
Denote $\sgn f(+\infty)=+1$ (resp. $-1$)
if $\sgn f(x)=+1$ (resp. $-1$) for sufficiently large $x$.
The meaning of $\sgn f(-\infty)$ is similar.
\begin{theorem}\label{Ffg}
Let $F,f,g$ be three real polynomials satisfying the following conditions.
\begin{enumerate}
\item [(a)]
$F(x)=a(x)f(x)+b(x)g(x)$,
where $a(x),b(x)$ are two real polynomials, such that
$\deg F=\deg f$ or $\deg f+1$.
\item [(b)]
$f,g\in\rz$ and $g\sep f$.
\item [(c)]
$F$ and $g$ have leading coefficients of the same sign.
\end{enumerate}
Suppose that $b(r)\le 0$ whenever $f(r)=0$.
Then $F\in\rz$ and $f\sep F$.
In particular, if $g\seps f$ and $b(r)<0$ whenever $f(r)=0$,
then $f\seps F$.
\end{theorem}
\begin{proof}
Without loss of generality,
we may assume that $f$ and $g$ have no common zeros,
which implies that $f$ has only simple zeros.
In other words, it suffices to consider the case $g\seps f$.
We may also assume that $F,g$ are standard
and that $\deg a(x)\le 1,\deg b(x)\le 2$.
If $b(x)\equiv 0$, then the result is trivial.
So let $b(x)\not\equiv 0$.

First consider the case $b(r)<0$ whenever $f(r)=0$.
Let $\deg f=n$ and $r_n<\cdots<r_1$ be all zeros of $f$.
Then $\sgn g(r_k)=(-1)^{k-1}$ since $g$ is standard and $g\seps f$.
It follows that $\sgn F(r_k)=(-1)^k$.
Since $F$ is standard,
we have $\sgn F(+\infty)=+1$
and $\sgn F(-\infty)=(-1)^{n+1}$ provided $\deg F=n+1$.
By Weierstrass Intermediate Value Theorem,
$F(x)$ has one zero in each of $n$ intervals
$(r_n,r_{n-1}),\ldots,(r_2,r_1),(r_1,+\infty)$
and has an additional zero in the interval $(-\infty,r_n)$ provided $\deg F=n+1$.
Hence $F$ has $\deg F$ real zeros and $f\seps F$.

For the general case,
define $b_j(x)=b(x)-{1}/{j}$ and $F_j(x)=a(x)f(x)+b_j(x)g(x)$.
Let $j$ be sufficiently large.
Then $b_j(r_k)<0$ for each zero $r_k$ of $f$
since the number of zeros of $b(x)$ is finite,
and so $F_j\in\rz$ and $f\seps F_j$ by the above discussion.
It is clear that $\deg F_j=\deg F$ since $F_j=F-g/j$.
Let $r_1>r_2>r_3>\cdots$ and $t_1^{(j)}>t_2^{(j)}>t_3^{(j)}>\cdots$
be all zeros of $f$ and $F_j$ respectively.
Then $t_1^{(j)}>r_1>t_2^{(j)}>r_2>t_3^{(j)}>r_3>\cdots$.
Note that the zeros of a polynomial are continuous functions of the coefficients of the polynomial
(see \cite{Coo08} for instance).
In particular the limit of a sequence of $\rz$ polynomials is still a $\rz$ polynomial.
Hence $F\in\rz$.
Moreover, let $t_1\ge t_2\ge t_3\ge\cdots$ be all zeros of $F$.
Then
$t_1\ge r_1\ge t_2\ge r_2\ge t_3\ge r_3\ge\cdots$
by continuity.
Thus $f\sep F$ and the proof is complete.
\end{proof}

Generally speaking, the conditions (a), (b) and (c) in
Theorem~\ref{Ffg} can be satisfied naturally. It remains to examine
the sign of $b(x)$ for the zeros of $f$. Sometimes this task is
trivial, for example, when $b(x)=-(b_0+b_1x)^2$. On the other hand,
if coefficients of $f$ are nonnegative (resp. alternating in sign),
then it suffices to consider the sign of $b(x)$ for $x\le 0$ (resp.
$x\ge 0$). As an example, we give a short and simple proof of
Haglund~\cite[Lemma 3.6]{Hag00}. We also refer the reader to Theorem
2.4.2--2.4.6 in Brenti~\cite{Bre89} for which Theorem \ref{Ffg} can
be used to give unified proofs.
\begin{corollary}[{\cite[Lemma 3.6]{Hag00}}]
Let $f$ and $g$ be two real polynomials with positive leading coefficients
$\alpha$ and $\beta$ respectively.
Suppose that the following conditions are satisfied.
\begin{enumerate}
\item [(a)]
$f,g\in\rz$ and $g$ interlaces $f$.
\item [(b)]
$F(x)=(ax+b)f(x)+x(x+d)g(x)$
where $a,b,d\in\mathbb{R}$ with $d\ge 0,d\ge b/a$,
and either $a>0$ or $a<-\beta/\alpha$.
\item [(c)]
All zeros of $f$ are nonpositive if $a>0$
and nonnegative if $a<-\beta/\alpha$.
\end{enumerate}
Then $F\in\rz$.
In addition,
if each zero $r$ of $f$ satisfies $-d\le r\le 0$,
then $f$ interlaces $F$.
\end{corollary}
\begin{proof}
Suppose that $a>0$.
Then $f\in\pf$, and so $g\in\pf$.
It follows from Corollary \ref{WY-cor} that $F\in\rz$ since $a\cdot d\ge b\cdot 1$.

Now suppose that $a<-\beta/\alpha$.
Then $f$ and $g$ have coefficients alternating in sign
and $F$ has negative leading coefficient.
Let $\deg f=n$. Then $\deg g=n-1$ and $\deg F=n+1$.
Define $f_1(x)=(-1)^nf(-x), g_1(x)=(-1)^{n-1}g(-x)$
and $F_1(x)=(-1)^{n}F(-x)$.
Then $f_1,g_1\in\pf$ and $g_1$ interlaces $f_1$.
Note that $F_1(x)=(-ax+b)f_1(x)+x(-x+d)g_1(x)$ and $(-a)\cdot d\ge b\cdot (-1)$.
Hence $F_1\in\rz$ by Corollary \ref{WY-cor},
and so $F\in\rz$.

Finally, if $-d\le r\le 0$ whenever $f(r)=0$,
then $f,g\in\pf$ and $r(r+d)\le 0$.
Thus $f$ interlaces $F$ by Theorem \ref{Ffg}.
\end{proof}

The following theorem is a generalization of Theorem~\ref{Ffg}
and gives a solution to Problem~\ref{prob-Ffgk}.
It can be proved by the same technique used in the proof of Theorem~\ref{Ffg}.
So we omit its proof for brevity.
\begin{theorem}\label{Ffgk}
Let $F,f,g_1,\ldots,g_k$ be real polynomials satisfying the following conditions.
\begin{enumerate}
\item [(a)]
$F(x)=a(x)f(x)+b_1(x)g_1(x)+\cdots+b_k(x)g_k(x)$,
where $a(x),b_1(x),\ldots,b_k(x)$ are real polynomials,
such that $deg F=\deg f$ or $\deg f+1$.
\item [(b)]
$f,g_j\in\rz$ and $g_j\sep f$ for each $j$.
\item [(c)]
$F$ and $g_1,\ldots,g_k$ have leading coefficients of the same sign.
\end{enumerate}
Suppose that $b_j(r)\le 0$
for each $j$ and each zero $r$ of $f$.
Then $F\in\rz$ and $f\sep F$.
In particular, if for each zero $r$ of $f$,
there is an index $j$ such that $g_j\seps f$ and $b_j(r)<0$,
then $f\seps F$.
\end{theorem}
\begin{corollary}\label{PP'P}
Let $\{P_n(x)\}$ be a sequence of standard polynomials and satisfy the recurrence relation
\begin{equation*}
P_{n}(x)=a_n(x)P_{n-1}(x)+b_n(x)P'_{n-1}(x)+c_n(x)P_{n-2}(x),
\end{equation*}
where $a_n(x),b_n(x),c_n(x)$ are real polynomials
such that $\deg P_n=\deg P_{n-1}$ or $\deg P_{n-1}+1$.
Suppose that for each $n$,
coefficients of $P_n(x)$ are nonnegative (resp. alternating in sign).
If $b_n(x)\le 0$ and $c_n(x)\le 0$ whenever $x\le 0$ (resp. $x\ge 0$),
then $\{P_n(x)\}$ forms a generalized Sturm sequence.
In particular,
if for each $n$, $\deg P_n=n$ and
either $b_n(x)<0$ or $c_n(x)<0$ whenever $x\le 0$ (resp. $x\ge 0$),
then $\{P_n(x)\}$ forms a Sturm sequence.
\end{corollary}

Corollary~\ref{PP'P} provides a unified approach to
the reality of zeros of certain well-known polynomials,
including the orthogonal polynomials $p_n(x)$ satisfying (\ref{op-rr}),
the Stirling polynomials $S_n(x)$ of the second kind satisfying (\ref{snx-rr})
and the Eulerian polynomials $A_n(x)$ satisfying (\ref{anq-rr}).
We will give some more applications of Corollary~\ref{PP'P} in the next section.
\begin{lemma}\label{Gfg}
Let $G(x)=c(x)f(x)+d(x)g(x)$
where $G,f,g$ are standard
and $c(x),d(x)$ are real polynomials.
Suppose that $f,g\in\rz$ and $g\seps f$.
Then the following hold.
\begin{enumerate}
\item [(i)]
If $\deg G\le\deg g+1$ and $c(s)>0$ whenever $g(s)=0$,
then $G\in\rz$ and $g\seps G$.
\item [(ii)]
If $\deg G\le\deg f$ and $d(r)>0$ whenever $f(r)=0$,
then $G\in\rz$ and $G\seps f$.
\end{enumerate}
The statements also hold if all instances of $\seps$ and $>$ are replaced by $\sep$ and $\ge$ respectively.
\end{lemma}
\begin{proof}
(i)\quad
Let $\deg g=m$ and $s_{m}<s_{m-1}<\cdots<s_2<s_1$ be all zeros of $g$.
Then $\sgn G(s_k)=(-1)^k$ since $c(s_k)>0$.
Also, $\sgn G(+\infty)=+1$,
and $\sgn G(-\infty)=(-1)^{m+1}$ provided $\deg G=m+1$.
Hence $G$ has one zero in each of $m$ intervals
$(s_m,s_{m-1})$, \ldots, $(s_2,s_1)$, $(s_1,+\infty)$
and has an additional zero in the interval $(-\infty,s_m)$ provided $\deg G=m+1$.
Thus $G\in\rz$ and $g\seps G$.

(ii)\quad
Let $\deg f=n$ and $r_n<r_{n-1}<\cdots<r_2<r_1$ be all zeros of $f$.
Then $\sgn G(r_k)=(-1)^{k-1}$, and $\sgn G(-\infty)=(-1)^n$ provided $\deg G=n$.
It follows that $G$ has one zero in each of $n-1$ intervals
$(r_n,r_{n-1})$, \ldots, $(r_2,r_1)$
and has an additional zero in the interval $(-\infty,r_n)$ provided $\deg G=n$.
Thus $G\in\rz$ and $G\seps f$.

The remaining part of the lemma follows by a continuity argument.
\end{proof}

The following theorem gives a solution to Problem~\ref{prob-FGfg}.
\begin{theorem}\label{FGfg}
Let $f,g,F,G$ be four standard real polynomials satisfying the following conditions.
\begin{enumerate}
\item [(a)]
$F(x)=a(x)f(x)+b(x)g(x)$ and $G(x)=c(x)f(x)+d(x)g(x)$
where $a(x)$, $b(x)$, $c(x)$, $d(x)$ are real polynomials
such that $\deg F=\deg G$ or $\deg G+1$.
\item [(b)]
$f,g\in\rz$ and $g\sep f$.
\item [(c)]
$\Delta(x):=a(x)d(x)-b(x)c(x)\ge 0$ whenever $G(x)=0$.
\end{enumerate}
Suppose that either $c(x)$ is a positive constant and $\deg G\le\deg g+1$
or $d(x)$ is a positive constant and $\deg G\le\deg f$.
Then $F,G\in\rz$ and $G\sep F$.
In particular, if $g\seps f$ and $\Delta(x)>0$ whenever $G(x)=0$, then $G\seps F$.
\end{theorem}
\begin{proof}
Suppose that $c(x)$ is a positive constant.
Then by Lemma~\ref{Gfg}~(i),
$g\sep f$ (resp. $g\seps f$) implies that $G\in\rz$ and $g\sep G$ (resp. $g\seps G$).
By Condition (a), we have $cF=aG+(bc-ad)g$.
If $\Delta(x)\ge 0$ (resp. $\Delta(x)>0$) whenever $G(x)=0$,
then by Theorem~\ref{Ffg}, $F\in\rz$ and $G\sep F$ (resp. $G\seps F$).

Suppose that $d(x)$ is a positive constant.
Then by Lemma~\ref{Gfg}~(ii),
$g\sep f$ (resp. $g\seps f$) implies that $G\in\rz$ and $G\sep f$ (resp. $G\seps f$).
By Condition (a), we have $dF=(ad-bc)f+bG$.
If $\Delta(x)\ge 0$ (resp. $\Delta(x)>0$) whenever $G(x)=0$,
then by Lemma~\ref{Gfg}~(i), $F\in\rz$ and $G\sep F$ (resp. $G\seps F$).
\end{proof}

When both $c$ and $d$ are constants,
Theorem~\ref{FGfg} is particularly interesting and useful as we shall see.
\begin{corollary}\label{cfdg}
Let $f,g,F,G$ be standard real polynomials and satisfy the following conditions.
\begin{enumerate}
\item [(a)]
$F(x)=a(x)f(x)+b(x)g(x)$ and $G(x)=cf(x)+dg(x)$
where $a(x),b(x)\in\mathbb{R}[x]$ and $c,d\in\mathbb{R}$.
\item [(b)]
$\deg F=\deg G$ or $\deg G+1$.
\item [(c)]
$f,g\in\rz$ and $g\sep f$.
\end{enumerate}
Suppose that $da(x)\ge cb(x)$ whenever $G(x)=0$.
Then $F,G\in\rz$ and $G\sep F$.
In particular,
if $g\seps f$ and $da(x)>cb(x)$ whenever $G(x)=0$,
then $G\seps F$.
\end{corollary}

It is well known that if $g\sep f$,
then $cf+dg\in\rz$ for arbitrary real numbers $c$ and $d$,
which is an immediate consequence of Theorem~\ref{WY-thm}
and can also be obtained by setting $a(x)\equiv c$ and $b(x)\equiv d$ in Corollary~\ref{cfdg}.
Using Corollary~\ref{cfdg} we can obtain more precise results.
The following are special cases of Corollary~\ref{cfdg}
when both $a$ and $b$ are constants.
Some of them have appeared in the literature,
for example, Wagner~\cite{Wag91,Wag92}.
\begin{corollary}\label{abcd}
Let $a,b,c,d\ge 0$.
Suppose that $f,g\in\rz$ are standard and $g\sep f$.
Then the following statements hold.
\begin{enumerate}
\item [(i)]
If $ad\ge bc$, then $cf+dg\sep af+bg$.
In particular, $g\sep af+bg\sep f$.
\item [(ii)]
If $af-bg$ is standard,
then $cf+dg\sep af-bg$.
In particular, $f,g\sep af-bg$.
\item [(iii)]
If $-af+bg$ is standard, then $-af+bg\sep cf+dg$,
and in particular, $-af+bg\sep f,g$.
\end{enumerate}
Similar results hold when $g\seps f$.
\end{corollary}

Using Corollary~\ref{cfdg} we can give a short and simple proof of Theorem~\ref{WY-thm}.
In fact, we can prove the following result more precise than Theorem~\ref{WY-thm}.
\begin{corollary}
Let $f,g\in\rz$ have leading coefficients of the same sign and $g\sep f$.
Define $F=(ax+b)f+(cx+d)g$ and $G=af+cg$ where $a,b,c,d\in\mathbb{R}$.
If $ad\le bc$,
then $F,G\in\rz$ and $G\sep F$.
\end{corollary}
\begin{proof}
We have $F=xG+H$ where $H=bf+dg$.
Clearly, $G,H\in\rz$.
If $G\equiv 0$, then the statement is trivial,
so let $G\not\equiv 0$.
Without loss of generality,
we may assume that both $f$ and $g$ are monic.
We may also assume that $G$ is standard
(otherwise replaced $G$ by $-G$).
Note that $(ax+b)c-(cx+d)a=bc-ad\ge 0$.
To prove the result by means of Corollary~\ref{cfdg},
it suffices to prove that $F$ is standard.

Actually,
if $a=0$ or $\deg G=\deg f$,
then $F$ is obviously standard.
Assume now that $a\neq 0$ and $\deg G<\deg f$.
Then $\deg f=\deg g$ and $a=-c$.
Let $f=\prod_{i=1}^{n}(x-r_i)$ and $g=\prod_{i=1}^{n}(x-s_i)$.
Then
$$G=af+cg=c\left(\sum_{i=1}^{n}r_i-\sum_{i=1}^{n}s_i\right)x^{n-1}+\cdots.$$
Note that $g$ alternates $f$ and $f\not\equiv g$ (otherwise $G\equiv 0$).
Hence $\sum_{i=1}^{n}s_i<\sum_{i=1}^{n}r_i$,
which yields that $c>0$ since $G$ is standard.
It follows that $b+d\ge 0$ from $ad\le bc$.
Thus $H=bf+dg$ is standard, and so is $F=xG+H$.
The proof is complete.
\end{proof}
\section{Applications and related topics}
\hspace*{\parindent}
In this section we apply the results obtained in the previous section
to derive several known facts and to solve certain new problems
in a unified manner.
\subsection{Orthogonal polynomials}
\hspace*{\parindent}
The classical orthogonal polynomials satisfy a three-term recurrence relation
\begin{equation*}
p_n(x)=(a_nx+b_n)p_{n-1}(x)-c_np_{n-2}(x)
\end{equation*}
with $p_0(x)=1$ and $p_1(x)=a_1x+b_1$,
where $a_n,c_n>0$ and $b_n\in\mathbb{R}$ for all $n\ge 1$
(see \cite[Theorem 3.2.1]{Sze75}).
By Theorem \ref{Ffg} or Corollary \ref{PP'P},
the orthogonal polynomials $p_n(x)$ form a Sturm sequence.
This is a standard result in the theory of orthogonal polynomials.
The following are some classical orthogonal polynomials
(see Szeg\"o~\cite{Sze75} for details).
\begin{itemize}
\item (Tchebyshev)\quad
$T_{n+1}(x)=2xT_n(x)-T_{n-1}(x),\quad T_1(x)=x$ or $T_1(x)=2x$.
\item (Hermite)\quad
$H_{n+1}(x)=2xH_n(x)-2nH_{n-1}(x),\quad H_1(x)=2x$.
\item (Laguerre)\quad
$(n+1)L_{n+1}(x)=(2n+1-x)L_n(x)-nL_{n-1}(x),\quad L_1(x)=1-x$.
\item (Legendre)\quad
$(n+1)P_{n+1}(x)=(2n+1)xP_n(x)-(n-1)P_{n-1}(x),\quad P_1(x)=x$.
\item (Gegenbauer)\quad For $\la>-1/2$,
$$(n+1)C_{n+1}^{\la}(x)=2(n+\la)xC_{n}^{\la}(x)-(n+2\la-1)C_{n-1}^{\la}(x),
\quad C_1^{\la}(x)=2\la x.$$
\item (Jacobi)\quad For $\alpha,\beta>-1$,
\begin{eqnarray*}
& & 2n(n+\alpha+\beta)(2n+\alpha+\beta-2)P_n^{(\alpha,\beta)}(x)\\
&=& (2n+\alpha+\beta-1)[(2n+\alpha+\beta)(2n+\alpha+\beta-2)x+\alpha^2-\beta^2]P_{n-1}^{(\alpha,\beta)}(x)\\
& & -2(n+\alpha-1)(n+\beta-1)(2n+\alpha+\beta)P_{n-2}^{(\alpha,\beta)}(x),\quad n\ge 2,
\end{eqnarray*}
with $P_1^{(\alpha,\beta)}(x)=\frac{1}{2}(\alpha+\beta+2)x+\frac{1}{2}(\alpha-\beta)$.

For $\alpha=\beta=0$, $P_n^{(0,0)}(x)$ reduces to the Legendre polynomials.
The Gegenbauer polynomials and the Tchebyshev polynomials
can also be viewed as special cases of the Jacobi polynomials.
\end{itemize}
Note that the leading coefficients of the Laguerre polynomials $L_n(x)$ have the sign $(-1)^n$.
Set ${\overline L}_n(x)=L_n(-x)$.
Then ${\overline L}_n(x)$ are standard and satisfy
$$(n+1){\overline L}_{n+1}(x)=(2n+1+x){\overline L}_n(x)-n{\overline L}_{n-1}(x).$$
Thus $\{{\overline L}_n(x)\}$ forms a Sturm sequence,
and so does $\{L_n(x)\}$.
\subsection{Matching polynomials}
\hspace*{\parindent}
Let $G$ be a graph with $n$ vertices
and $p(G,k)$ the number of matchings of size $k$,
i.e., the number of sets of $k$ edges of $G$,
no two edges having a common vertex.
Set $p(G,0)=1$ for convenience.
The matching polynomial
$M(G,x)=\sum_k(-1)^{k}p(G,k)x^{n-2k}$
counts the matchings in the graph $G$.
Clearly, for any $v\in V(G)$,
$$p(G,k)=p(G-\{v\},k)+\sum_{u\sim v}p(G-\{v,u\},k-1),$$
where the first term in the sum counts $k$-matchings which do not use $v$
and the second one counts $k$-matching which do use $v$.
This leads to a recurrence relation
\begin{equation*}
M(G,x)=xM(G-\{v\},x)-\sum_{u\sim v}M(G-\{v,u\},x).
\end{equation*}
From Theorem \ref{Ffgk} and by induction on the number of vertices of graphs,
it is easy to see that $M(G,x)$ has only real zeros and $M(G-\{v\},x)$ interlaces $M(G,x)$ for any $v\in V(G)$,
a well-known result
(see Godsil and Gutman~\cite{GG81} for instance).

The matching polynomials,
formally introduced by Farrell~\cite{Far79} in 1979,
have occurred not only in the combinatorial literature
but also in various branches of physics and chemistry
(see \cite{GG81} for a brief survey).
For example, to understand the behavior of a monomer-dimer system in statistical physics,
Heilmann and Lieb~\cite{HL72} in 1972
introduced the partition function $Q(G,x)$ of a monomer-dimer system,
which is essentially the matching polynomial of the graph $G$
with the edge weight $W$.
If $W\equiv 1$,
then $Q(G,x)$ is precisely the ordinary matching polynomial of the graph $G$.
It is easy to yield the recurrence relation
$$Q(G,x)=xQ(G-\{v\},x)-\sum_{u\sim v}W(u,v)Q(G-\{v,u\},x).$$
Heilmann and Lieb~\cite[Theorem 4.2]{HL72} showed that
if $G$ is a simple graph with nonnegative edge weights,
then $Q(G,x)$ has only real zeros and
$Q(G-\{v\},x)$ interlaces $Q(G,x)$ for any $v\in V(G)$.
This result is now clear from the point of view of Theorem~\ref{Ffgk}.

The classical orthogonal polynomials are closely related to the matching polynomials.
For example,
the Tchebyshev polynomials of two kinds are the matching polynomials of paths and cycles respectively,
the Hermite polynomials and the Laguerre polynomials
are the matching polynomials of completes graphs and complete bipartite graphs respectively.

Another class of polynomials closely related to the matching polynomials is the rook polynomials.
The original rook polynomial defined in Riordan's book~\cite{Rio58}
counts the number of ways of arranging nonattacking rooks on a board
(a board is a finite subset of $\mathbb{N}\times\mathbb{N}$).
Goldman, Joichi and White~\cite{GJW76} conjectured that such rook polynomials have only real zeros.
Wilf extended the concept of rook polynomial to the case of arbitrary matrix instead of a board
(a board is identified with a $(0,1)$-matrix)
and Nijenhuis showed that the extended rook polynomials of nonnegative matrices have only real zeros
by establishing a result similar to Theorem~\ref{Ffgk} (see \cite{Nij76}).
Bender observed that the (extended) rook polynomial
is the same as the matching polynomial of a (weighted) bipartite graph,
and Nijenhuis's result and Goldman-Joichi-White's conjecture
are therefore implied by the result of Heilmann and Lieb.
\subsection{Brenti's derangement polynomials}
\hspace*{\parindent}
Let $S_n$ denote the symmetric group on $n$-elements
$[n]=\{1,2,\ldots,n\}$.
Let $\pi$ be a permutation in $S_n$.
An element $i\in [n]$ is called an excedance of $\pi$ if $\pi(i)>i$.
Denote by $\exc(\pi)$ the number of excedances of $\pi$.
The permutation $\pi$ is called a derangement if $\pi(i)\ne i$ for all $i\in [n]$.
Let $D_n$ denote the set of all derangements of $S_n$.
Define the derangement polynomial
$$d_n(q)=\sum_{\pi\in D_n}q^{\exc(\pi)}.$$
For example,
$d_0(q)=1,d_1(q)=0,d_2(q)=q,d_3(q)=q+q^2$. 
Since $d_n(1)=|D_n|$ we may consider
$d_n(q)$ as a $q$-analogue of the derangement numbers.

The exponential generating function of $d_n(q)$ can be written as
\begin{equation}\label{dnq-ef}
\sum_{n\ge 0}d_n(q)\dfrac{t^n}{n!}=\dfrac{1}{1-\sum_{n\ge 2}(q+q^2+\cdots+q^{n-1})t^n/n!}
\end{equation}
(see \cite{Bre90,Ros68} for instance).
Gessel~\cite{Ges91} gave a direct proof of (\ref{dnq-ef}) with $q=1$
based on a factorization of some so-called D-permutations.
Kim and Zeng~\cite{KZ01} gave a decomposition of derangements
which interprets (\ref{dnq-ef}) directly.
Using (\ref{dnq-ef})
Brenti~\cite{Bre90} showed that
the polynomial $d_n(q)$ is symmetric and unimodal for $n\ge 1$.
He further proposed the following.
\begin{conjecture}[{\cite[Conjecture]{Bre90}\label{dp-conj}}]
The polynomial $d_n(q)$ has only real zeros for $n\ge 1$.
\end{conjecture}

The derangement polynomials are closely related to the classical Eulerian polynomials
which are defined by
$$A_n(q)=\sum_{\pi\in S_n}q^{\exc(\pi)+1}$$
for $n\ge 1$ and $A_0(q)=1$
(see Stanley~\cite[Proposition 1.3.12]{Sta97I} for instance).
Therefore,
$$A_n(q)/q
=\sum_{s\subseteq [n]}\sum_{\pi\in D_{|s|}}q^{\exc(\pi)}
=\sum_{k=0}^{n}\binom{n}{k}d_k(q).$$
From the binomial inversion formula it follows that
$$d_n(q)=\sum_{k=0}^{n}(-1)^{n-k}\binom{n}{k}A_k(q)/q.$$
Recall that the Eulerian polynomials satisfy the recurrence relation
$$A_n(q)=nqA_{n-1}(q)+q(1-q)A'_{n-1}(q).$$
Hence $\{d_{n}(q)\}$ satisfies the recurrence relation
\begin{equation}\label{dnq-rr}
d_{n}(q)=(n-1)qd_{n-1}(q)+q(1-q)d'_{n-1}(q)+(n-1)qd_{n-2}(q).
\end{equation}
Zhang~\cite{Zhang-Hefei,Zhang-Xi'an} verified and generalized Conjecture~\ref{dp-conj} as follows:

``Let $f_n(q)$ be a polynomial of degree $n$ with nonnegative real coefficients
and satisfy
\begin{itemize}
\item [(a)]
$f_n(q)=a_nqf_{n-1}(q)+b_nq(1+c_nq)f'_{n-1}(q)+d_nqf_{n-2}(q)$,
where $a_n>0,b_n>0,d_n\ge 0$ and $n\ge 2$.
\item [(b)]
For $n\ge 1$, zero is a simple root of $f_n(q)$.
\item [(c)]
$f_0(q)=e, f_1(q)=e_1q$ and $f_2(q)$ has two real roots,
where $e\ge 0,e_1\ge 0$.
\end{itemize}
Then the polynomial $f_n(q)$ has $n$ distinct real roots,
separated by the roots of $f_{n-1}(q)$.''

However, Zhang's result does not hold in general
and his proof is valid only for $\{d_n(q)\}$.
The condition $c_n\le 0$ is necessary in Zhang's result,
even in the case $d_n=0$.
For example,
let $f_0(q)=1,f_1(q)=q$,
$f_2(q)=qf_1(q)+q(q+2)f'_1(q)=2q(q+1)$ and
$f_3(q)=qf_2(q)+q(2q+1)f'_3(q)=2q(5q^2+5q+1)$.
Then zeros of $f_2(q)$ are $-1$ and $0$
and those of $f_3(q)$ are $(-5\pm\sqrt{5})/10$ and $0$.
Clearly, the latter cannot be separated by the former.

Now from the point of view of Corollary~\ref{PP'P},
Conjecture~\ref{dp-conj} is obviously true
and furthermore,
the derangement polynomials $d_n(q)$ form a generalized Sturm sequence.
\begin{remark}
The referee pointed out that Conjecture~\ref{dp-conj} has been proved by Canfield (unpublished).
\end{remark}
\subsection{Narayana polynomials}
\hspace*{\parindent}
The Narayana numbers $N(n,k)=\frac{1}{n}\binom{n}{k}\binom{n}{k-1}$
and the Narayana polynomials $N_n(q)=\sum_{k=1}^{n}N(n,k)q^k$
have many combinatorial interpretations and fascinating properties
(see \cite{Sul99,Sul00,Sul02} for instance).
For example,
the Narayana numbers $N(n,k)$ can be viewed as a refinement of the famous Catalan numbers $C_n$
since $C_n=N_n(1)$
(see R\'emy~\cite{Rem85} for a combinatorial proof of this and
Stanley~\cite{Sta97II,Sta-cn} for various combinatorial interpretations of the Catalan numbers).
Also, $N_n(2)$ are the Schr\"oder numbers
(see Foata and Zeilberger \cite{FZ97} for a combinatorial proof
and Stanley \cite{Sta97AMM} for an interesting history of the Schr\"oder numbers).
In fact, the original Narayana polynomials ${\overline N}_n(q)$,
introduced by Bonin, Shapiro and Simion~\cite{BSS93},
are defined as the $q$-analog of the Schr\"oder numbers.
It is shown in \cite{BSS93} that ${\overline N}_n(q)$ has unimodal coefficients and
has $q=-1$ as one zero for $n\ge 1$ and $q=-2$ as one zero for even $n\ge 2$.
It is known that ${\overline N}_n(q)=N_n(1+q)$ and
\begin{equation}\label{nnx-rr}
(n+1)N_n(q)=(2n-1)(1+q)N_{n-1}(q)-(n-2)(1-q)^2N_{n-2}(q)
\end{equation}
with $N_1(q)=q$ and $N_2(q)=q(1+q)$
(see \cite{Sul02} for instance).
As pointed out by Stanley (see B\'ona~\cite{Bon02}),
Theorem 5.3.1 in Brenti~\cite{Bre89} implies that $N_n(q)$ have only real zeros.
Another proof for this result was recently found by Br\"and\'en~\cite{Bra0303149}
by expressing $N_n(q)$ in terms of the Jacobi polynomials:
$$N_{n}(q)=\frac{1}{n+1}(1-q)^{n}P_n^{(1,1)}\left(\frac{1+q}{1-q}\right).$$

More generally,
Sulanke~\cite{Sul99} defined polynomial sequences $\{p_{\alpha,n}(q)\}_{n\ge 2}$,
for $\alpha=0,1,2$,
in terms of various ``diagonal thickness'' parameters on parallelogram polyominoes.
These polynomials satisfy the recurrence relation
$$(n+1-\alpha)p_{\alpha,n+1}(q)=
(2n-1-\alpha)(1+q)p_{\alpha,n}(q)-(n-2)(1-q)^2p_{\alpha,n-1}(q)$$
with $p_{\alpha,2}(q)=q$ and $p_{\alpha,3}(q)=q(1+q)$,
and are generalizations of some well-known combinatorial sequences.
In particular, $p_{0,n}(q)=N_n(q)$.

Using Theorem~\ref{Ffg} or Corollary~\ref{PP'P},
we can give a more direct and natural approach to
the reality of zeros of $N_n(q)$ as well as $p_{\alpha,n}(q)$.
\begin{proposition}
Let $\{P_n(x)\}$ be a sequence of polynomials with nonnegative coefficients
and $\deg P_n=\deg P_{n-1}+1$.
Suppose that
$$P_n(x)=(a_nx+b_n)P_{n-1}(x)-(c_nx+d_n)^2P_{n-2}(x)$$
where $a_n,b_n,c_n,d_n\in\mathbb{R}$.
Then $\{P_n(x)\}$ forms a Sturm sequence.
\end{proposition}

It is well known that
the Catalan numbers (the Narayana numbers) count
the number of $1$-stack sortable $n$-permutations (with $k$ descents).
The problem of stack sorting was introduced by Knuth~\cite{Knu73} in 1960's
and many variations have been considered since then
(see B\'ona~\cite{Bon03} for a survey).
Let $W_t(n,k)$ be the number of $t$-stack sortable $n$-permutations with $k$ descents
and $W_{n,t}(q)=\sum_{k=0}^{n-1}W_t(n,k)q^k$ the associated generating function.
In \cite{Bon02}, B\'ona showed that $\{W_t(n,k)\}_{0\le k\le n-1}$ is a unimodal sequence for fixed $n$ and $t$
and further conjectured that $W_{n,t}(q)$ have only real zeros for $n\ge 2$.
Note that $W_{n,1}(q)=N_n(q)$, the Narayana polynomial,
and $W_{n,n-1}(q)=A_n(q)$, the classical Eulerian polynomial.
Hence the conjecture is true for $t=1,n-1$.
Br\"and\'en~\cite{BraTAMS} has recently verified the conjecture for $t=2,n-2$.
But the conjecture remains open in the general case.
\subsection{Compositions of multisets and Dowling lattices}
\hspace*{\parindent}
Let $\bn=(n_1,n_2,\ldots)$ be the multiset consisting of $n_i$ copies of the $i$th type element.
Denote by $\mo(\bn,k)$ the number of compositions of $\bn$ into exactly $k$ parts.
Then
\begin{equation}\label{onk-rr}
(n_j+1)\mo(\bn+e_j,k)=k\mo(\bn,k-1)+(n_j+k)\mo(\bn,k),
\end{equation}
where $\bn+e_j$ denotes the multiset obtained from $\bn$
by adjoining one (additional) copy of the $j$th type element
(see Riordan~\cite[p.96]{Rio58}).
Let $f_{\bn}(x)=\sum_{k\ge 0}\mo(\bn,k)x^k$ be the associated generating function.
Using (\ref{onk-rr}) Simion~\cite{Sim84} deduced the recurrence relation
\begin{equation}\label{fnejx-rr}
(n_j+1)f_{\bn+e_j}(x)=(x+n_j)f_{\bn}(x)+x(x+1)f'_{\bn}(x).
\end{equation}
By means of appropriate transformation to (\ref{fnejx-rr}),
Simion~\cite[Theorem 1]{Sim84} can show that
the polynomial $f_{\bn}(x)$ has all zeros in the interval $[-1,0]$,
and furthermore, $f_{\bn}(x)$ and $f_{\bn+e_j}(x)$ have interlaced zeros.
This result is now clear from the point of view of Theorem~\ref{Ffg}.

In particular,
if $\bn=(1,1,\ldots,1)$, then $\mo(\bn,k)=k!S(n,k)$
where $S(n,k)$ is the Stirling number of the second kind.
Thus the polynomial $F_n(x)=\sum_{k=1}^{n}k!S(n,k)x^k$ has only real zeros,
which can also be followed from the formula
$F_n(x)=\frac{x^{n+1}}{x+1}A_n(\frac{x+1}{x})$
where $A_n(x)$ is the Eulerian polynomial
(see \cite{Ben97} for instance).
The polynomial $F_n(x)$
was first studied by Tanny~\cite{Tan75}.
Benoumhani~\cite{Ben97} gave a generalization of $F_n(x)$
replacing the Stiling numbers by the Whitney numbers of the Dowling lattices.
The Dowling lattice $Q_n(G)$ is a geometric lattice of rank $n$ over a finite group $G$ of order $m$
and have many remarkable properties
(see \cite{Ben96,Ben97,Ben99,Dow73} for instance).
When $m=1$, that is, $G$ is the trivial group,
$Q_n(G)$ is isomorphic to
the lattice $\Pi_{n+1}$ of partitions of an $(n+1)$-element set.
So the Dowling lattices can be viewed as group-theoretic analogs of the partition lattices.
Let $W_m(n,k)$ be the $k$th Whitney numbers of the second kind of $Q_n(G)$.
Denote $D_m(n;x)=\sum_{k=0}^{n}W_m(n,k)x^k$ and $F_n(m;x)=\sum_{k=0}^{n}k!W_m(n,k)x^k$.
(Dowling gave a combinatorial interpretation for the coefficients $k!W_m(n,k)$, see \cite{Ben97}.)
Then
$$D_m(n;x)=(x+1)D_m(n-1;x)+mxD'_m(n-1;x)$$
(see \cite{Ben99}) and
$$F_n(m;x)=(x+1)F_{n-1}(m;x)+x(x+m)F'_{n-1}(m;x)$$
(see \cite{Ben97}).
Benoumhani showed that
both $D_m(n;x)$ and $F_n(m;x)$ have only real zeros for $n\ge 1$
(see \cite[Theorem 2]{Ben99} and \cite[Theorem 6]{Ben97} respectively).
These results can also be followed from Corollary~\ref{WY-cor}.
As a consequence, $W_m(n,k)$ is unimodal and log-concave in $k$.
This gives supports to a long-standing conjecture that
the Whitney numbers of the second kind of any finite geometric lattice
are unimodal or even log-concave
(see \cite[Conjecture 3]{Sta89}).
When $m=1$,
we have $Q_n(G)\cong \Pi_{n+1}$ and $W_m(n,k)=S(n,k)$.
Again we obtain that the polynomials
$S_n(x)=\sum_{k=0}^{n}S(n,k)x^k$ and $F_n(x)=\sum_{k=0}^{n}k!S(n,k)x^k$
have only real zeros for $n\ge 1$.
\subsection{Eulerian polynomials of Coxeter groups}
\hspace*{\parindent}
Given a finite Coxeter group $W$,
define the Eulerian polynomials of $W$ by
$$P(W,x)=\sum_{\pi\in W}x^{d_W(\pi)},$$
where $d_W(\pi)$ is the number of $W$-descents of $\pi$.
We refer the reader to Bj\"orner and Brenti \cite{BB05} for relevant definitions.
Brenti proposed the following conjecture.
\begin{conjecture}[{\cite[Conjecture 5.2]{Bre94EuJC}}]\label{pwx-conj}
For every finite Coxeter group $W$,
the polynomial $P(W;x)$ has only real zeros.
\end{conjecture}

For Coxeter groups of type $A_n$,
it is known that $P(A_n,x)=A_n(x)/x$, the shifted Eulerian polynomial.
The classical Eulerian polynomial $A_n(x)$ has only real zeros,
so does $P(A_n,x)$.
Since $\{A_n(x)\}$ satisfies
$$A_n(x)=nxA_{n-1}(x)+x(1-x)A'_{n-1}(x),\quad A_0(x)=1,$$
it immediately yields that $\{P(A_n,x)\}$ satisfies
$$P(A_n,x)=(nx+1)P(A_{n-1},x)+x(1-x)P'(A_{n-1},x),\quad P(A_0,x)=1.$$

In \cite{FS70}, Foata and Sch\"utzenberger introduced
a $q$-analog of the classical Eulerian polynomials defined by
$$A_n(x;q)=\sum_{\pi\in\msn}x^{\exc(\pi)+1}q^{c(\pi)},$$
where $\exc(\pi)$ and $c(\pi)$ denote the numbers of excedances and cycles in $\pi$ respectively.
It is clear that $A_n(x;1)=A_n(x)$ is precisely the classical Eulerian polynomial.
Brenti showed that $q$-Eulerian polynomials satisfy the recurrence relation
$$A_{n}(x;q)=(nx+q-1)A_{n-1}(x;q)+x(1-x)\frac{\partial}{\partial x}A_{n-1}(x;q),$$
with $A_0(x;q)=x$
(\cite[Proposition 7.2]{Bre00}).
He showed also that
$A_n(x;q)$ has only real nonnegative simple zeros when $q$ is a positive rational number
(\cite[Theorem 7.5]{Bre00}).

For Coxeter groups of type $B_n$,
Brenti~\cite{Bre94EuJC} defined a $q$-analogues of $P(B_n,x)$,
which reduces to $A_n(x)$ when $q=0$ and to $P(B_n,q)$ when $q=1$, by
$$B_n(x;q)=\sum_{\pi\in B_n}q^{N(\pi)}x^{d_B(\pi)}$$
where $N(\pi)=|\{i\in [n]: \pi(i)<0\}|$.
He showed that $\{B_n(x;q)\}$ satisfies the recurrence relation
$$B_n(x;q)=\{1+[(1+q)n-1]x\}B_{n-1}(x;q)+(1+q)x(1-x)\frac{\partial}{\partial x}B'_{n-1}(x;q),$$
with $B_0(x;q)=1$ (\cite[Theorem 3.4 (i)]{Bre94EuJC})
and that all $B_n(x;q)$ have only real zeros for $q\ge 0$
(\cite[Corollary 3.7]{Bre94EuJC}).
In particular $P(B_n,x)$ has only real zeros.

By the classification of finite irreducible Coxeter groups,
it suffices to decide whether the conjecture holds for Coxeter groups of type $D_n$
to settle Conjecture \ref{pwx-conj}
(see Brenti~\cite{Bre94EuJC} for further information).

Using Theorem~\ref{Ffg} or Corollary~\ref{PP'P},
we can give a unified interpretation of
the reality of zeros of $A_n(x),P(A_n,x),A_n(x;q)$ and $B_n(x;q)$.
\begin{proposition}
Let $\{P_n(x)\}$ be a sequence of polynomials with nonnegative coefficients
and $\deg P_n=\deg P_{n-1}+1$.
Suppose that
$$P_n(x)=(a_nx+b_n)P_{n-1}(x)+x(c_nx+d_n)P'_{n-1}(x)$$
where $a_n,b_n\in\mathbb{R}$ and $c_n\le 0,d_n\ge 0$.
Then $\{P_n(x)\}$ forms a generalized Sturm sequence.
\end{proposition}
\subsection{Genus polynomials of graphs}
\hspace*{\parindent}
Given a finite graph $G$ and a nonnegative integer $k$,
let $\gamma(G,k)$ denote the number of distinct embeddings of the graph $G$
into an oriented surface of genus $k$.
We refer the reader to \cite{GT87} for the basic terminology of graph embeddings.
The genus polynomial is defined in \cite{GF87} as $GP(G,x)=\sum_{k\ge 0}\gamma(G,k)x^k$.
Gross, Robbins and Tucker~\cite{GRT89}
showed that the genus distribution of the bouquet is log-concave
and conjectured that the genus distribution of every graph is log-concave
(i.e., $\gamma(G,k)$ is log-concave in $k$).
Stahl~\cite[Conjecture 6.4]{Stahl97} further conjectured
that the genus polynomial of every graph has only real zeros.
He also verified the conjecture for several infinite families of graphs
by establishing the recurrence relations of the associated genus polynomials.
These results follow from our results in the previous section.
In particular,
Stahl considered the $H$-linear family of graphs
obtained by consistently amalgamating additional copies of a graph $H$.
For such a family $\{G_n\}$,
there is a square matrix $M$ and a vector $v$ with entries in $\mathbb{Z}[x]$
such that the genus polynomial of $G_n$ is the first entry of $M^nv$
(\cite[Proposition 5.2]{Stahl97}).
The following are genus generating matrices and initial vectors of certain linear families
given by Stahl~\cite{Stahl97}.
\begin{example}\label{ex-lf}
\mbox{}
\begin{enumerate}
\item [(i)]
For the cobblestone paths,
$M_1=\begin{pmatrix} 4 & 2\\ 6x & 0\end{pmatrix}$
and $v_1=\begin{pmatrix} 1\\ x \end{pmatrix}$.
\item [(ii)]
For the ladders,
$M_2=\begin{pmatrix} 0 & 4\\ 2x & 2\end{pmatrix}$
$v_2=\begin{pmatrix} 1\\ 1 \end{pmatrix}$.
\item [(iii)]
For the double ladders,
$M_3=6\begin{pmatrix} 3x & 3\\ 2x & 1+3x\end{pmatrix}$
and $v_3=2\begin{pmatrix} 2\\ 1+x \end{pmatrix}$.
\item [(iv)]
For the diamonds,
$M_4=4\begin{pmatrix} 2+3x & 1\\ 4x & 2x\end{pmatrix}$
and $v_4=2\begin{pmatrix} 1+x\\ 2x \end{pmatrix}$.
\item [(v)]
For the triple ladders,
$M_5=\begin{pmatrix} 192x & 96+288x\\ 72+192x^2 & 24+288x\end{pmatrix}$
and $v_5=\begin{pmatrix} 18+18x\\ 6+30x \end{pmatrix}$.
\item [(vi)]
For the $K_4$-linear graphs,
$M_6=\begin{pmatrix} 8+68x & 4+16x\\ 32x+48x^2 & 16x\end{pmatrix}$
and $v_6=\begin{pmatrix} 2+14x\\ 8x+8x^2 \end{pmatrix}$.
\item [(vii)]
For the $W_4$-linear graphs,
$$M_7=4\begin{pmatrix} 2+65x+54x^2 & 1+22x\\ 16x+104x^2 & 8x+16x^2\end{pmatrix}
\text{ and } v_7=\begin{pmatrix} 2+58x+36x^2\\ 16+80x \end{pmatrix}.$$
\item [(viii)]
For the triangular prisms,
$$M_8=\begin{pmatrix} 0 & 162x & 54\\ 24x^2 & 72x & 12+108x\\ 11x^2 & 15x+117x^2 & 1+72x\end{pmatrix}
\text{ and }v_8=\begin{pmatrix} 8\\ 4+4x\\ 1+7x \end{pmatrix}.$$
\end{enumerate}
\end{example}

Stahl verified his conjecture for the linear families
associated with Example~\ref{ex-lf} (i), (ii) and (iv),
i.e., the first entry of $M^kv$ has only real zeros for $k\ge 1$.
He left the remaining as a conjecture and asked some more general questions as follows.
\begin{conjecture}[{\cite[Conjecture 6.9]{Stahl97}}]\label{stahl-conj}
The zeros of the genus polynomials of the graphs listed in Example \ref{ex-lf} are real and negative.
\end{conjecture}
\begin{question}[{\cite[Question 6.10 and 6.11]{Stahl97}}]\label{stahl-q}
Let $M(x)$ be a square matrix whose entries are real polynomials.
\begin{enumerate}
\item [(i)]
Under what conditions,
if $(f(x),g(x))$ is a pair of polynomials whose zeros interlace,
do the zeros of the two components of the vector $(f(x),g(x))M(x)$ interlace?
\item [(ii)]
Under what conditions,
are the zeros of each of the entries of $M^k(x)$
all real for $k=1,2,\ldots$?
\end{enumerate}
\end{question}

Stahl showed that the matrix $M_4=4\begin{pmatrix} 2+3x & 1\\ 4x & 2x\end{pmatrix}$
in Example \ref{ex-lf} has both properties in Question \ref{stahl-q},
but the matrix $M=\begin{pmatrix} 3x & 3\\ 2x & 3x+1000\end{pmatrix}$ has neither.
In what follows,
we apply Theorem~\ref{FGfg} to give an answer to Question~\ref{stahl-q}
and to verify Conjecture~\ref{stahl-conj} for the graphs in Example~\ref{ex-lf} (i)-(vi)
in a unified approach.
\begin{definition}
Let $M=\begin{pmatrix} a(x) & c(x)\\ b(x) & d(x)\end{pmatrix}$,
where $a,b,c,d$ are polynomials with nonnegative coefficients.
We say that the polynomial matrix $M$ is {\it nice} if
\begin{enumerate}
\item [(a)]
$\deg a,\deg d\le 1, \deg b\le 2$
and $c$ is a positive constant.
\item [(b)]
$\det(M)\ge 0$ for $x\le 0$.
\end{enumerate}
\end{definition}
\begin{proposition}\label{stahl-answer}
Let $M$ be a nice matrix.
Suppose that $f,g\in\pf$ and $g\sep f$.
Then the following hold.
\begin{enumerate}
\item [(i)]
If $(F,G)=(f,g)M$, then $F,G\in\pf$ and $G\sep F$.
\item [(ii)]
If $(G_1,F_1)^T=M(g,f)^T$, then $F_1,G_1\in\pf$ and $G_1\sep F_1$.
\item [(iii)]
Each entry of $M^k$ has only real zeros for $k=1,2,\ldots$.
\end{enumerate}
\end{proposition}
\begin{proof}
We have $F=af+bg$ and $G=cf+dg$.
We may assume that $F\not\equiv 0$ and $G\not\equiv 0$.
To prove (i) by means of Theorem~\ref{FGfg},
it suffices to prove that $\deg F=\deg G$ or $\deg G+1$.
This is obvious if $\deg f=\deg g+1$.
So let $\deg f=\deg g=n$.
We distinguish two cases.
Suppose first that $\deg G=n+1$.
Then $\deg d=1$.
Since $ad-bc\ge 0$ whenever $x\le 0$,
we have $\deg a\ge 1$ or $\deg b\ge 1$,
Hence $\deg F\ge n+1$.
On the other hand,
it is clear that $\deg F\le n+2$ since $\deg a\le 1$ and $\deg b\le 2$.
Suppose now that $\deg G=n$.
Then $d$ is a constant.
Again by the assumption $ad-bc\ge 0$ whenever $x\le 0$,
we have $\deg b\le 1$.
It follows that $n\le\deg F\le n+1$.
Thus (i) follows from Theorem~\ref{FGfg}.

Similarly, (ii) follows from Theorem~\ref{FGfg}
since $F_1=df+bg$ and $G_1=cf+ag$.

Finally, we apply (i) to prove (iii).
For $k=1,2,\ldots$, let
$$M^k=\begin{pmatrix} f_k^{(1)} & g_k^{(1)}\\ f_k^{(2)} & g_k^{(2)}\end{pmatrix}.$$
Then $(f_{1}^{(1)},g_{1}^{(1)})=(a,c), (f_{1}^{(2)},g_{1}^{(2)})=(b,d)$
and $(f_{k+1}^{(i)},g_{k+1}^{(i)})=(f_{k}^{(i)},g_{k}^{(i)})M$ for $i=1,2$.
By the assumption $ad-bc\ge 0$ whenever $x\le 0$
it follows that $d(r)=0$ implies $b(r)\le 0$.
This means that $b\in\pf$ and $d\sep b$.
Also, $c\sep a$ since $\deg c=0$ and $\deg a\le 1$.
Thus (iii) follows from (i) by induction on $k$.
\end{proof}
\begin{proposition}
The zeros of the genus polynomials of the graphs listed in Example~\ref{ex-lf}~(i)-(vi)
are real and negative.
\end{proposition}
\begin{proof}
We need to verify that for $i=1,2,\ldots,6$,
the first entry of the column vector $M_i^kv_i$ has only real zeros for each $k$.
By Proposition \ref{stahl-answer} (ii) and by induction on $k$,
it suffices to prove that each $M_i$ is the product of certain nice matrices
and $v_i^{(1)}\sep v_i^{(2)}$ where $v_i=(v_i^{(1)},v_i^{(2)})^T$.

For $i=1,2,3,4$,
it is easy to verify that each matrix $M_i$ is nice.
For $i=5,6$, the matrix $M_i$ can be decomposed into the product of certain nice matrices:
\begin{eqnarray*}
M_5&=&24\begin{pmatrix} 4+12x & 8\\ 1+12x & 3+8x\end{pmatrix}\begin{pmatrix} 0 & 1\\ x & 0\end{pmatrix},\\
M_6&=&4\begin{pmatrix} 0 & 1\\ x & 0\end{pmatrix}\begin{pmatrix} 8+12x & 4\\ 2+17x & 1+4x\end{pmatrix}
\end{eqnarray*}

It is also easy to verify that $v_i^{(1)}\sep v_i^{(2)}$ for $i=1,2,\ldots,6$.
Thus the proof is complete.
\end{proof}
\begin{remark}
The genus polynomials of the $W_4$-linear graphs in Example~\ref{ex-lf}~(vii)
may have nonreal zeros.
For example, let $u_7=M_7v_7$. Then the first entry of $u_7$ is
$$u_7^{(1)}=8(10+339x+2855x^2+2736x^3+972x^4).$$
Using Mathematica,
we can obtain the approximations of four zeros of $u_7^{(1)}$:
$$x^{(1)}_{1,2}=-1.34194\pm i0.88376,\quad
x^{(1)}_3=-0.0828403,\quad
x^{(1)}_4=-0.0481022.$$
This gives a counterexample to Conjecture~\ref{stahl-conj}.
Stahl's conjecture about the reality of zeros of the genus polynomials
is therefore false in general.
But it is possible that the genus distribution of each graph is log-concave.
\end{remark}

We end this paper by proposing the following.
\begin{problem}
Characterize all real polynomial matrices that can be decomposed into the product of finite nice matrices
and find an algorithm of decomposition for such matrices.
\end{problem}
\section*{Acknowledgements}
\hspace*{\parindent}
This paper was completed when the second author's stay at the
Institute of Mathematics, Academia Sinica, Taipei. He thanks Prof.
K.-W. Lih and Prof. Y.-N. Yeh for their kind assistance in the
preparation of this paper.

The authors thank the anonymous referee for his/her careful reading
and valuable suggestions.

\end{document}